\documentclass[12pt]{article}%
\usepackage{graphicx}
\usepackage[intlimits]{amsmath}
\usepackage{latexsym}
\usepackage{amsfonts}
\usepackage{amssymb}%
\makeatletter
\newfont{\footsc}{cmcsc10 at 8truept}
\newfont{\footbf}{cmbx10 at 8truept}
\newfont{\footrm}{cmr10 at 10truept}
\newtheorem{theorem}{Theorem}

\newtheorem{lemma}[theorem]{Lemma}

\newenvironment{proof}[1][Proof]{\noindent{\textbf {#1}  }}  {\hfill$\Box$\bigskip}

\begin{document}

\title{Bounds on graph eigenvalues I }
\author{Vladimir Nikiforov\\Department of Mathematical Sciences, University of Memphis, \\Memphis TN 38152, USA}
\maketitle

\begin{abstract}
We improve some recent results on graph eigenvalues. In particular, we prove
that if $G$ is a graph of order $n\geq2,$ maximum degree $\Delta,$ and girth
at least $5,$ then
\[
\mu\left(  G\right)  \leq\min\left\{  \Delta,\sqrt{n-1}\right\}  ,
\]
where $\mu\left(  G\right)  $ is the largest eigenvalue of the adjacency
matrix of $G$.

Also, if $G$ is a graph of order $n\geq2$ with dominating number
$\gamma\left(  G\right)  =\gamma,$ then%
\begin{align*}
\lambda_{2}\left(  G\right)   &  \leq\left\{
\begin{array}
[c]{cc}%
n & \text{if }\gamma=1\\
n-\gamma & \text{if }\gamma\geq2,
\end{array}
\right. \\
\lambda_{n}\left(  G\right)   &  \geq\left\lceil n/\gamma\right\rceil ,
\end{align*}
where $0=\lambda_{1}\left(  G\right)  \leq\lambda_{2}\left(  G\right)
\leq...\leq\lambda_{n}\left(  G\right)  $ are the eigenvalues of the Laplacian
of $G.$

We also determine all cases of equality in the above inequalities.



\end{abstract}

\section{Introduction}

Our notation is standard (e.g., see \cite{Bol98} and \cite{CDS80}); in
particular, we write $G\left(  n\right)  $ for a graph of order $n$ and
$G\left(  n,m\right)  $ for a graph of order $n$ and $m$ edges. Given a vertex
$u\in V\left(  G\right)  ,$ we write $\Gamma\left(  u\right)  $ for the set of
neighbors of $u$ and set $d\left(  u\right)  =\left\vert \Gamma\left(
u\right)  \right\vert .$ If $X,Y$ are two disjoint subsets of $V\left(
G\right)  ,$ we denote by $e\left(  X,Y\right)  $ the number of $X-Y$ edges.
Given a graph $G$ of order $n,$ we write $\mu\left(  G\right)  =\mu_{1}\left(
G\right)  \geq...\geq\mu_{n}\left(  G\right)  $ for the eigenvalues of its
adjacency matrix and $0=\lambda_{1}\left(  G\right)  \leq...\leq\lambda
_{n}\left(  G\right)  =\lambda\left(  G\right)  $ for the eigenvalues of its Laplacian.

This note is motivated by some recent papers on graph eigenvalues. Liu, Lu,
and Tian \cite{LLT05b} proved that if $G=G\left(  n\right)  $ is a connected
graph of girth at least $5$ and maximum degree $\Delta,$ then%
\begin{equation}
\mu\left(  G\right)  \leq\frac{-1+\sqrt{4n+4\Delta-3}}{2}; \label{main1}%
\end{equation}
equality holds if and only if $G=C_{5}.$

Observe that equality holds in (\ref{main1}) also for $K_{2}$ and all $\Delta
$-regular Moore graphs of diameter $2$. Hoffman and Singleton \cite{HoSi68}
proved that $r$-regular Moore graphs of diameter $2$ exist for $r=2,3,7$ and
possibly $57.$

A stronger theorem follows from a result in \cite{FMS93}.

\begin{theorem}
\label{th1}Let $G=G\left(  n\right)  $ be a graph of maximum degree $\Delta$
and girth at least $5.$ Then%
\begin{equation}
\mu\left(  G\right)  \leq\min\left\{  \Delta,\sqrt{n-1}\right\}  .
\label{main2}%
\end{equation}
Equality holds if and only if one of the following conditions holds:

(i) $G=K_{1,n-1};$

(ii) $G$ is a $\Delta$-regular Moore graph of diameter $2$;

(iii) $G=G_{1}\cup G_{2},$ where $G_{1}$ is $\Delta$-regular, $\Delta\left(
G_{2}\right)  \leq\Delta,$ and the girth of both $G_{1}$ and $G_{2}$ is at
least $5.$
\end{theorem}

Note that the right-hand side of (\ref{main2}) never exceeds the right-hand
side of (\ref{main1}).

Given a graph $G,$ a set $X\subset V\left(  G\right)  $ is called
\emph{dominating,} if $\Gamma\left(  u\right)  \cap X\neq\varnothing$ for
every $u\in V\left(  G\right)  \backslash X.$ The number $\gamma\left(
G\right)  =\min\left\{  \left\vert X\right\vert :X\text{ is a dominating
set}\right\}  $ is called the \emph{dominating number} of $G.$

Liu, Lu, and Tian \cite{LLT05a} proved\ that if $n\geq2$ and $G=G\left(
n\right)  $ is a connected graph with $\gamma\left(  G\right)  =\gamma,$ then%
\[
\lambda_{2}\left(  G\right)  \leq n-\gamma+\frac{n-\gamma^{2}}{n-\gamma}.
\]

For $n\geq\gamma^{2}$ this bound is implied by the following theorem.

\begin{theorem}
\label{th2}Let $n\geq2$ and $G=G\left(  n\right)  $ be a graph with
$\gamma\left(  G\right)  =\gamma.$ Then%
\begin{equation}
\lambda_{2}\left(  G\right)  \leq\left\{
\begin{array}
[c]{cc}%
n & \text{if }\gamma=1\\
n-\gamma & \text{if }\gamma\geq2.
\end{array}
\right.  \label{main5}%
\end{equation}
If $\gamma=1,$ equality holds if and only if $G=K_{n}$. If $\gamma=2,$
equality holds if and only if $G$ is the complement of a perfect matching. If
$\gamma>2,$ (\ref{main5}) is always a strict inequality.
\end{theorem}

Another result of Liu, Lu, and Tian \cite{LLT05a} states that if $n\geq2$ and
$G=G\left(  n\right)  $ is a connected graph with $\gamma\left(  G\right)
=\gamma,$ then%
\begin{equation}
\lambda\left(  G\right)  \geq n/\gamma; \label{main4}%
\end{equation}
equality holds if and only if $K_{1,n-1}\subset G.$

Inequality (\ref{main4}) follows immediately from a known result stated in
Lemma 4 of the same paper - Mohar \cite{Moh92} proved that for every set
$X\subset V=V\left(  G\right)  ,$ the inequality $\lambda\left(  G\right)
\left\vert X\right\vert \left\vert V\backslash X\right\vert \geq ne\left(
X,V\backslash X\right)  $ holds. Hence, if $X$ is a dominating set with
$\left\vert X\right\vert =\gamma,$ then $e\left(  X,V\backslash X\right)
\geq\left\vert V\backslash X\right\vert =n-\gamma,$ and (\ref{main4}) follows.

In fact, a subtler theorem holds.

\begin{theorem}
\label{th3}Let $n\geq2$ and $G=G\left(  n\right)  $ be a graph with
$\gamma\left(  G\right)  =\gamma>0.$ Then%
\begin{equation}
\lambda\left(  G\right)  \geq\left\lceil n/\gamma\right\rceil . \label{main6}%
\end{equation}
Equality holds if and only if $G=G_{1}\cup G_{2},$ where $G_{1}$ and $G_{2}$
satisfy the following conditions:

(i) $\left\vert G_{1}\right\vert =\left\lceil n/\gamma\right\rceil $ and
$\gamma\left(  G_{1}\right)  =1$;

(ii) $\gamma\left(  G_{2}\right)  =\gamma-1$ and $\lambda\left(  G_{2}\right)
\leq\left\lceil n/\gamma\right\rceil .$
\end{theorem}

Note that the above results of Liu, Lu, and Tian are stated for connected
graphs only, illustrating a tendency in some papers on graph eigenvalues to
stipulate connectedness \textit{apriori} - see, e.g., \cite{DaKu04},
\cite{Pan02}, \cite{SHW02}, \cite{ShWu04}, \cite{YLT04}, \cite{Zha04}. If not
truly necessary, such stipulation sends a wrong message. For example, Hong
\cite{Hon93} stated his famous inequality
\[
\mu\left(  G\right)  \leq\sqrt{2e\left(  G\right)  -v\left(  G\right)  +1}%
\]
for connected graphs, although his proof works for graphs with minimum degree
at least $1.$ This result has been reproduced verbatim countless times
challenging the readers to complete the picture on their own. Confinement to
connected graphs simplifies the study of cases of equality, but important
points might be missed. As an illustration, recall the result of Hong, Shu,
and Fang \cite{HSF01}: if $G=G\left(  n,m\right)  $ is a connected graph with
$\delta\left(  G\right)  =\delta,$ then
\begin{equation}
\mu\left(  G\right)  \leq\frac{\delta-1+\sqrt{8m-4\delta n+\left(
\delta+1\right)  ^{2}}}{2}, \label{mainin}%
\end{equation}
with equality holding if and only if every vertex of $G$ has degree $\delta$
or $n-1.$

Inequality (\ref{mainin}) has been proved independently by Nikiforov
\cite{Nik02} for disconnected graphs as well; however, as shown in
\cite{Nik02} and \cite{ZhCh05}, there are nonobvious disconnected graphs for
which equality holds in (\ref{mainin}).

Finally, observe that (\ref{mainin}) implies a result of Cao \cite{Cao98},
recently reproved for connected graphs by Das and Kumar \cite{DaKu04}: if
$G=G\left(  n,m\right)  $ is a graph with $\delta\left(  G\right)  =\delta
\geq1$ and $\Delta\left(  G\right)  =\Delta,$ then
\begin{equation}
\mu\left(  G\right)  \leq\sqrt{2m-\left(  n-1\right)  \delta+\left(
\delta-1\right)  \Delta}. \label{maincao}%
\end{equation}
In fact, (\ref{maincao}) follows from (\ref{mainin}) by%
\[
\mu^{2}\left(  G\right)  \leq2m-\left(  n-1\right)  \delta+\left(
\delta-1\right)  \mu\left(  G\right)  \leq2m-\left(  n-1\right)
\delta+\left(  \delta-1\right)  \Delta.
\]

\section{Proofs}

We shall need the following result of Grone and Merris \cite{GrMe94}.

\begin{lemma}
\label{leGM}If $G$ is a graph with $e\left(  G\right)  >0$, then
$\lambda\left(  G\right)  \geq\Delta\left(  G\right)  +1;$ if $G$ is
connected, then equality holds if and only if $\Delta\left(  G\right)
=\left\vert G\right\vert -1$.
\end{lemma}

\begin{proof}
[\textbf{Proof of Theorem \ref{th1}}]Since $\mu\left(  G\right)  \leq
\Delta\left(  G\right)  ,$ to prove (\ref{main2}), all we need is to show that
$\mu\left(  G\right)  \leq\sqrt{n-1}.$ We follow here the argument of Favaron,
Mah\'{e}o, and Sacl\'{e} \cite{FMS93}, p. 203. For every $u\in V\left(
G\right)  $ set
\[
w\left(  u\right)  =\sum_{v\in\Gamma\left(  u\right)  }d\left(  v\right)  .
\]
As shown in \cite{FMS93}, p. 203, $\mu^{2}\left(  G\right)  \leq\max_{u\in
V\left(  G\right)  }w\left(  u\right)  ;$ if $G$ is connected, equality holds
if and only if $G$ is regular or bipartite semiregular graph.

We shall prove that $w\left(  u\right)  \leq n-1$ for every $u\in V\left(
G\right)  .$ Indeed, let $u\in V\left(  G\right)  ;$ for every two distinct
vertices $v,w\in\Gamma\left(  u\right)  ,$ in view of $C_{3}\nsubseteq G$ and
$C_{4}\nsubseteq G,$ we see that $e\left(  \Gamma\left(  u\right)  \right)
=0$\ and $\Gamma\left(  v\right)  \cap\Gamma\left(  w\right)  =\left\{
u\right\}  .$ Hence,
\begin{equation}
w\left(  u\right)  =\sum_{v\in\Gamma\left(  u\right)  }d\left(  v\right)
=e\left(  \Gamma\left(  u\right)  ,V\left(  G\right)  \backslash\Gamma\left(
u\right)  \right)  \leq d\left(  u\right)  +n-d\left(  u\right)  -1=n-1,
\label{in2}%
\end{equation}
completing the proof of (\ref{main2}). Note that if equality holds in
(\ref{in2}), then $\cup_{v\in\Gamma\left(  u\right)  }\Gamma\left(  v\right)
=V\left(  G\right)  ;$ hence, the distance of any vertex $v\in V\left(
G\right)  $ to $u$ is at most $2$.

Let us determine when equality holds in (\ref{main2}). If any of the
conditions \emph{(i)-(iii) }holds, clearly (\ref{main2}) is an equality.
Suppose equality holds in (\ref{main2}). If $\mu\left(  G\right)  =\Delta,$
then $G$ contains a $\Delta$-regular component, say $G_{1}.$ Writing $G_{2}$
for the union of the remaining components of $G,$ we see that \emph{(iii)
}holds, completing the proof in this case.

Now let $\mu\left(  G\right)  =\sqrt{n-1}$; hence, equality holds in
(\ref{in2}) for some vertex $u\in V\left(  G\right)  ,$ implying that $G$ is
connected. According to the aforementioned result of Favaron, Mah\'{e}o, and
Sacl\'{e}, equality in (\ref{in2}) holds for every vertex $u\in V\left(
G\right)  ,$ and $G$ is either $\sqrt{n-1}$-regular or bipartite semiregular.
Clearly, $diam$ $G=2,$ so if $G$ is $\sqrt{n-1}$-regular, then it is a Moore
graph and \emph{(ii)} holds.

Finally, let $G$ be a bipartite graph. Then the distance between any two
vertices belonging to different parts of $G$ is odd; since $diam$ $G=2,$ it
follows that $G$ is a complete bipartite graph, and so, $G=K_{1,n-1},$
completing the proof.
\end{proof}

\bigskip

\begin{proof}
[\textbf{Proof of Theorem \ref{th2}}]Let $V=V\left(  G\right)  $ and $X$ be a
dominating set with $\left\vert X\right\vert =\gamma.$ For the sake of
completeness we shall reprove the known inequality $\delta\left(  G\right)
\leq n-\gamma.$ Indeed, select $u\in X.$ If $\Gamma\left(  u\right)  \cap
X=\varnothing,$ then $\Gamma\left(  u\right)  \subset V\backslash X,$ and so
$\delta\left(  G\right)  \leq d\left(  u\right)  \leq n-\gamma.$ Now assume
that $\Gamma\left(  u\right)  \cap X\neq\varnothing.$ Then there exists
$v\in\left(  V\backslash X\right)  \cap\Gamma\left(  u\right)  $ such that $v$
is not joined to any $w\in X\backslash\left\{  u\right\}  ,$ otherwise
$X\backslash\left\{  u\right\}  $ would be dominating, contradicting that $X$
is minimal. Hence,
\[
\delta\left(  G\right)  \leq d\left(  v\right)  \leq\left(  n-1\right)
-\left\vert X\backslash\left\{  u\right\}  \right\vert =n-\gamma,
\]
as claimed.

If $G=K_{n},$ we have $\gamma=1$ and $\lambda_{2}\left(  K_{n}\right)  =n,$
completing the proof. Assume that $G\neq K_{n};$ hence, $e\left(  \overline
{G}\right)  \geq1.$ Applying Lemma \ref{leGM}, we have%
\begin{equation}
\lambda_{2}\left(  G\right)  =n-\lambda\left(  \overline{G}\right)  \leq
n-\Delta\left(  \overline{G}\right)  -1=n-\left(  n-1-\delta\left(  G\right)
\right)  -1\leq n-\gamma, \label{in3}%
\end{equation}
proving (\ref{main5}).

Let us determine when equality holds in (\ref{main5}). If $\gamma=1$ and
$G=K_{n},$ then $\lambda_{2}\left(  G\right)  =n,$ so (\ref{main5}) is an
equality. If $\gamma=2$ and $G=\overline{\left(  n/2\right)  K_{2}}$ then
$\lambda_{2}\left(  G\right)  =n-\lambda\left(  \overline{G}\right)  =n-2,$ so
(\ref{main5}) is an equality.

Suppose now that equality holds in (\ref{main5}). If $\gamma=1,$ from
$\lambda_{2}\left(  G\right)  =n-\lambda\left(  \overline{G}\right)  =n$ we
find that $e\left(  \overline{G}\right)  =0,$ and so $G=K_{n}.$ If $\gamma
\geq2,$ then we have equalities in (\ref{in3}), implying that $\delta\left(
G\right)  =n-\gamma$ and $\lambda\left(  \overline{G}\right)  =\Delta\left(
\overline{G}\right)  +1=\gamma.$ From Lemma \ref{leGM} we conclude that
$\overline{G}$ has a component $G_{1}$ such that $\Delta\left(  G_{1}\right)
=\gamma-1$ and $\left\vert G_{1}\right\vert =\gamma.$ Set $V_{1}=V\left(
G_{1}\right)  .$ Since $G_{1}$ is a component of $\overline{G},$ the pair
$\left(  V_{1},V\backslash V_{1}\right)  $ induces a complete bipartite graph
in $G$ and so $\gamma=2.$ We have $\lambda\left(  \overline{G}\right)
=n-\lambda_{2}\left(  G\right)  =2$ and so $\Delta\left(  \overline{G}\right)
=1.$ This implies that $\overline{G}$ is a perfect matching, as otherwise $G$
would have a dominating vertex, contradicting that $\gamma=2.$ This completes
the proof.\bigskip
\end{proof}

\begin{proof}
[\textbf{Proof of Theorem \ref{th3}}]Let $G$ be a graph with $\gamma\left(
G\right)  =\gamma;$ set $V=V\left(  G\right)  $ and let $G_{0}\subset G$ be an
edge-minimal subgraph of $G$ with $V\left(  G_{0}\right)  =V$ and
$\gamma\left(  G_{0}\right)  =\gamma.$ Clearly, $G_{0}$ is a union of $\gamma$
vertex-disjoint stars, and so $G_{0}$ contains a star of order at least
$\left\lceil n/\gamma\right\rceil .$ Therefore, $\lambda\left(  G\right)
\geq\lambda\left(  G_{0}\right)  \geq\left\lceil n/\gamma\right\rceil ,$
proving (\ref{main6}).

Let us determine when equality holds in (\ref{main6}). If $G=G_{1}\cup G_{2},$
where $G_{1}$ and $G_{2}$ satisfy conditions \emph{(i)} and \emph{(ii)} of
Theorem \ref{th3}, then clearly equality holds in (\ref{main6}). Let $G$ be a
graph such that equality holds in (\ref{main6}), and $G_{0}\subset G$ be an
edge-minimal subgraph with $V\left(  G_{0}\right)  =V$ and $\gamma\left(
G\right)  =\gamma;$ clearly, $G_{0}$ is a union of $\gamma$ vertex-disjoint
stars, whose centers form a dominating set of $G$. From
\[
\left\lceil n/\gamma\right\rceil =\lambda\left(  G\right)  \geq\lambda\left(
G_{0}\right)  \geq\left\lceil n/\gamma\right\rceil
\]
we conclude that $G_{0}$ contains a component $H$ that is star
$K_{1,\left\lceil n/\gamma\right\rceil -1}.$ To complete the proof, we have to
show that no edge of $G$ joins $H$ to another component $F$ of $G_{0}.$ If
there is such an edge, according to Lemma \ref{leGM}, the component
$G^{\prime}$ of $G$ containing both $H$ and $F$ must satisfy $\lambda\left(
G^{\prime}\right)  >\left\lceil n/\gamma\right\rceil ,$ a contradiction.
Hence, $H$ induces a component of $G$, say $G_{1}.$ We have $\left\vert
G_{1}\right\vert =\gamma+1$ and $\gamma\left(  G_{1}\right)  =1,$ so
\emph{(i)} holds. Setting $G_{2}$ for the union of the remaining components of
$G$, we see that $\gamma\left(  G_{2}\right)  =\gamma-1,$ since $G_{2}$ is
spanned by $\gamma-1$ stars. Observing that $\lambda\left(  G_{2}\right)
\leq\lambda\left(  G\right)  =\left\lceil n/\gamma\right\rceil ,$ condition
\emph{(ii)} follows, completing the proof.
\end{proof}

\textbf{Acknowledgement }Lihua Feng and the referee pointed out some
shortcomings in an earlier version of the note.

\end{document}